\expandafter\ifx\csname mthreemacsloaded\endcsname\relax\else \fi

\magnification1100
\input amstex


 \catcode`\@=11
 \let\wlog@ld\wlog
 \def\wlog#1{\relax}

 \newif\ifIN@
 \def\m@rker{\m@@rker}
 \def\IN@{\expandafter\INN@\expandafter}
 \long\def\INN@0#1@#2@{\long\def\NI@##1#1##2##3\ENDNI@
    {\ifx\m@rker##2\IN@false\else\IN@true\fi}%
     \expandafter\NI@#2@@#1\m@rker\ENDNI@}
  \newtoks\Initialtoks@  \newtoks\Terminaltoks@
  \def\SPLIT@{\expandafter\SPLITT@\expandafter}
  \def\SPLITT@0#1@#2@{\def\TTILPS@##1#1##2@{%
     \Initialtoks@{##1}\Terminaltoks@{##2}}\expandafter\TTILPS@#2@}
  \newtoks\Trimtoks@

 \def\ForeTrim@{\expandafter\ForeTrim@@\expandafter}
 \def\ForePrim@0 #1@{\Trimtoks@{#1}}
 \def\ForeTrim@@0#1@{\IN@0\m@rker. @\m@rker.#1@%
     \ifIN@\ForePrim@0#1@%
     \else\Trimtoks@\expandafter{#1}\fi}
 
  \def\Trim@0#1@{%
      \ForeTrim@0#1@%
      \IN@0 @\the\Trimtoks@ @%
        \ifIN@
             \SPLIT@0 @\the\Trimtoks@ @\Trimtoks@\Initialtoks@
             \IN@0\the\Terminaltoks@ @ @%
                 \ifIN@
                 \else \Trimtoks@ {FigNameWithSpace}%
                 \fi
        \fi
      }

  \font\titlebold=cmbx12 scaled 1200
  \font\twelvebold=cmbx12
  \font\tenbold=cmbx10
  \font\ninebold=cmbx9
  \font\sevenbold=cmbx7
  \font\fivebold=cmbx5

  \input amssym.def \input amssym
     \font\titlemsa=msam10 at 14.4pt
     \font\titlemsb=msbm10 at 14.4pt
     \font\titleeufm=eufm10 at 14.4pt
     \font\twelvemsa=msam10 scaled 1200
     \font\twelvemsb=msbm10 scaled 1200
     \font\twelveeufm=eufm10 scaled 1200
     \font\ninemsa=msam9
     \font\ninemsb=msbm9
     \font\nineeufm=eufm9

   \ifx\cyrfam\undefined
   \else
     \immediate\write16{}%
     \message{ !!! cyr fonts already defined. !!! }
     \message{ --- edit out superfluous font defs? }
   \fi
   \newfam\cyrfam
       \font\titlecyr=wncyr10 scaled 1440 
       \font\twelvecyr=wncyr10 scaled 1200
       \font\tencyr=wncyr10
       \font\ninecyr=wncyr9
       \font\sevencyr=wncyr7
       \font\sixcyr=wncyr6

   \newfam\eusmfam
       \font\titleeusm=eusm10 scaled 1440
       \font\twelveeusm=eusm10 scaled 1200
       \font\teneusm=eusm10
       \font\nineeusm=eusm9
       \font\seveneusm=eusm7
       
       \font\fiveeusm=eusm5

\let\Cal\cal

    \font\ninemrm=cmr9 
    \font\ninei=cmmi9
    \font\ninesy=cmsy9 
    \skewchar\ninei='177
    \skewchar\ninesy='60

  \font\twelvemrm=cmr10 at 12pt 
  \font\twelvei=cmmi10 at 12pt
  \font\twelvesy=cmsy10 at 12pt

  \font\titlemrm=cmr10 at 14.4pt 
  \font\titlei=cmmi10 at 14.4pt
  \font\titlesy=cmsy10 at 14.4pt


  \def\Smallfonts{\ninepoint}

  \def\Hfont{\titlepoint\bf}
  \def\Authorfont{\twelvepoint\it}
  \def\HHfont{\twelvepoint\bf}
  \def\HHHfont{\bf}
  \def\Bibfont{\tenbf}
  \def\Coordfont{\nineit }

  \def \thfont {\bf }
  \def \pffont {\it\itSpacing }
  \def \rkfont {\bf }
  \def \dffont {\bf }
  \def \egfont {\bf }

 \def\ninepoint{%
  \def\rm{\fam0\ninerm}%
    \textfont0=\ninemrm  \scriptfont0=\sevenrm  \scriptscriptfont0=\fiverm
    \textfont1=\ninei    \scriptfont1=\seveni   \scriptscriptfont1=\fivei
  \def\mit{\fam1\ninei}%
  \def\oldstyle{\fam1\ninei}%
    \textfont2=\ninesy   \scriptfont2=\sevensy  \scriptscriptfont2=\fivesy
    \textfont3=\tenex    \scriptfont3=\tenex    \scriptscriptfont3=\tenex
  \def\it{\fam\itfam\nineit}%
    \textfont\itfam=\nineit
  \def\bf{\ifmmode\fam\bffam\else\ninebf\fi}%
    \textfont\bffam=\ninebold 
    \scriptfont\bffam=\sevenbold 
    \scriptscriptfont\bffam=\fivebold%
  \def\msa{\fam\msafam\ninemsa}%
    \textfont\msafam=\ninemsa 
    \scriptfont\msafam=\sevenmsa
    \scriptscriptfont\msafam=\fivemsa%
  \def\msb{\fam\msbfam\ninemsb}%
    \textfont\msbfam=\ninemsb%
    \scriptfont\msbfam=\sevenmsb%
    \scriptscriptfont\msbfam=\fivemsb%
  \def\eufm{\fam\eufmfam\nineeufm}%
    \textfont\eufmfam=\nineeufm
    \scriptfont\eufmfam=\seveneufm
    \scriptscriptfont\eufmfam=\fiveeufm
   \def\eusm{\fam\eusmfam\nineeusm}%
     \textfont\eusmfam=\nineeusm
     \scriptfont\eusmfam=\seveneusm
     \scriptscriptfont\eusmfam=\fiveeusm
   \def\cyr{\fam\cyrfam\ninecyr}%
     \textfont\cyrfam=\ninecyr
     \scriptfont\cyrfam=\sevencyr
     \scriptscriptfont\cyrfam=\sixcyr
  \setbox\strutbox=\hbox{\vrule
      height7pt depth3pt width0pt}%
   \baselineskip=10.8pt\rm}

 \let\eightpoint\ninepoint 

 \def\tenpoint{%
  \def\rm{\fam0\tenrm}%
    \textfont0=\tenmrm \scriptfont0=\sevenrm \scriptscriptfont0=\fiverm%
  \def\mit{\fam1\teni}%
  \def\oldstyle{\fam1\teni}%
    \textfont1=\teni   \scriptfont1=\seveni  \scriptscriptfont1=\fivei%
    \textfont2=\tensy  \scriptfont2=\sevensy \scriptscriptfont2=\fivesy%
    \textfont3=\tenex  \scriptfont3=\tenex   \scriptscriptfont3=\tenex%
  \def\it{\fam\itfam\tenit}%
    \textfont\itfam=\tenit%
  \def\bf{\ifmmode\fam\bffam\else\tenbf\fi}%
    \textfont\bffam=\tenbold
    \scriptfont\bffam=\sevenbold%
    \scriptscriptfont\bffam=\fivebold%
  \def\msa{\fam\msafam\tenmsa}%
    \textfont\msafam=\tenmsa%
    \scriptfont\msafam=\sevenmsa%
    \scriptscriptfont\msafam=\fivemsa%
  \def\msb{\fam\msbfam\tenmsb}%
    \textfont\msbfam=\tenmsb%
    \scriptfont\msbfam=\sevenmsb%
    \scriptscriptfont\msbfam=\fivemsb%
  \def\eufm{\fam\eufmfam\teneufm}%
   \textfont\eufmfam=\teneufm
   \scriptfont\eufmfam=\seveneufm
   \scriptscriptfont\eufmfam=\fiveeufm
   \def\eusm{\fam\eusmfam\teneusm}%
    \textfont\eusmfam=\teneusm
    \scriptfont\eusmfam=\seveneusm
    \scriptscriptfont\eusmfam=\fiveeusm
   \def\cyr{\fam\cyrfam\tencyr}%
    \textfont\cyrfam=\tencyr
    \scriptfont\cyrfam=\sevencyr
    \scriptscriptfont\cyrfam=\sixcyr
  \setbox\strutbox=\hbox{\vrule %
      height8.5pt depth3.5ptwidth0pt}%
  \baselineskip=\StdBaselineskip\rm}

 \def\twelvepoint{%
  \def\rm{\fam0\twelverm}%
    \textfont0=\twelvemrm \scriptfont0=\tenmrm \scriptscriptfont0=\sevenrm
    \textfont1=\twelvei   \scriptfont1=\teni   \scriptscriptfont1=\seveni
  \def\mit{\fam1\twelvei}%
  \def\oldstyle{\fam1\twelvei}%
    \textfont2=\twelvesy  \scriptfont2=\tensy  \scriptscriptfont2=\sevensy
    \textfont3=\tenex  \scriptfont3=\tenex  \scriptscriptfont3=\tenex
  \def\it{\fam\itfam\twelveit}%
    \textfont\itfam=\twelveit
  \def\bf{\ifmmode\fam\bffam\else\twelvebf\fi}%
    \textfont\bffam=\twelvebold
    \scriptfont\bffam=\tenbold%
    \scriptscriptfont\bffam=\sevenbold%
  \def\msa{\fam\msafam\twelvemsa}%
    \textfont\msafam=\twelvemsa%
    \scriptfont\msafam=\tenmsa%
    \scriptscriptfont\msafam=\sevenmsa%
  \def\msb{\fam\msbfam\twelvemsb}%
    \textfont\msbfam=\twelvemsb%
    \scriptfont\msbfam=\tenmsb%
    \scriptscriptfont\msbfam=\sevenmsb%
  \def\eufm{\fam\eufmfam\twelveeufm}%
   \textfont\eufmfam=\twelveeufm
   \scriptfont\eufmfam=\teneufm
   \scriptscriptfont\eufmfam=\seveneufm
   \def\eusm{\fam\eusmfam\twelveeusm}%
    \textfont\eusmfam=\twelveeusm
    \scriptfont\eusmfam=\teneusm
    \scriptscriptfont\eusmfam=\seveneusm
   \def\cyr{\fam\cyrfam\tencyr}%
    \textfont\cyrfam=\twelvecyr
    \scriptfont\cyrfam=\tencyr
    \scriptscriptfont\cyrfam=\sevencyr
  \setbox\strutbox=\hbox{\vrule
      height10.2pt depth4.55pt width0pt}%
  \baselineskip=14pt\rm}

 \def\titlepoint{%
    \textfont0=\titlemrm \scriptfont0=\twelvemrm \scriptscriptfont0=\tenmrm
    \textfont1=\titlei   \scriptfont1=\twelvei   \scriptscriptfont1=\teni
  \def\mit{\fam1\titlei}%
  \def\oldstyle{\fam1\titlei}%
    \textfont2=\titlesy  \scriptfont2=\twelvesy  \scriptscriptfont2=\tensy
    \textfont3=\tenex
    \scriptfont3=\tenex
    \scriptscriptfont3=\tenex
  \def\it{\fam\itfam\titleit}%
    \textfont\itfam=\titleit
  \def\bf{\ifmmode\fam\bffam\else\titlebf\fi}%
    \textfont\bffam=\titlebold
    \scriptfont\bffam=\twelvebold%
    \scriptscriptfont\bffam=\tenbold%
  \def\msa{\fam\msafam\titlemsa}%
    \textfont\msafam=\titlemsa%
    \scriptfont\msafam=\twelvemsa%
    \scriptscriptfont\msafam=\tenmsa%
  \def\msb{\fam\msbfam\titlemsb}%
    \textfont\msbfam=\titlemsb%
    \scriptfont\msbfam=\twelvemsb%
    \scriptscriptfont\msbfam=\tenmsb%
  \def\eufm{\fam\eufmfam\titleeufm}%
    \textfont\eufmfam=\titleeufm
    \scriptfont\eufmfam=\twelveeufm
    \scriptscriptfont\eufmfam=\teneufm
   \def\eusm{\fam\eusmfam\titleeusm}%
     \textfont\eusmfam=\titleeusm
     \scriptfont\eusmfam=\twelveeusm
     \scriptscriptfont\eusmfam=\teneusm
   \def\cyr{\fam\cyrfam\tencyr}%
    \textfont\cyrfam=\titlecyr
    \scriptfont\cyrfam=\twelvecyr
    \scriptscriptfont\cyrfam=\tencyr
  \setbox\strutbox=\hbox{\vrule
      height12.3pt depth5.54pt width0pt}%
  \baselineskip=16pt\rm}

\newbox\AuthorBox\newbox\TitleBox
\newbox\TFLinebox
\newbox\FLinebox
\newbox\HLinebox
\def\SetTFLinebox#1{\setbox\TFLinebox=\hbox{#1}}
\def\SetFLinebox#1{\setbox\FLinebox=\hbox{#1}}
\def\SetHLinebox#1{\setbox\HLinebox=\hbox{#1}}

 \def\SetAuthorHead#1{%
     \setbox\AuthorBox=\hbox{\ninepoint \it 
           \ignorespaces\frenchspacing#1\unskip}}
 \def\SetTitleHead#1{%
     \setbox\TitleBox=\hbox{\ninepoint \it
           \ignorespaces\frenchspacing#1\unskip}}

  \def\itSpacing{\relax}
  \def\itSpacingOff{\relax}


 \def\Hrule{\hrule width0pt height0pt}

  \newskip\ProcSkip \ProcSkip 8pt plus2pt minus2pt

 \newskip\LastSkip
 \def\SaveLastSkip{\LastSkip\lastskip}
 \def\RestoreLastSkip{\vskip-\LastSkip\vskip\LastSkip}

 \def\NoindentAfter{\everypar={\setbox0=\lastbox\everypar={}}}

 \long\def\H#1\par#2\par{\notenumber=0 \titlepagetrue%
    {
    \baselineskip=20pt
    \parindent=0pt\parskip=0pt\frenchspacing
    \leftskip=0pt plus .2\hsize minus .3\hsize
    \rightskip=0pt plus .2\hsize minus .3\hsize
 \def\\{\unskip\break}%
    \pretolerance=10000 \Hfont #1\unskip\break
     \vskip7pt\Hrule
\hfill \Authorfont #2\hfill\hfill\unskip}
    \vskip48pt plus 4pt minus 4pt
    \par\NoindentAfter\rm}

 \long\def\Hi#1\par#2\par{\notenumber=0 \titlepagetrue%
    {  \baselineskip=0pt  \parindent=0pt\parskip=0pt\frenchspacing
    \leftskip=0pt plus .2\hsize minus .3\hsize
    \rightskip=0pt plus .2\hsize minus .3\hsize
}
    \rm}


 \newdimen\PageRemainder
  \def\SetPageRemainder{
     \PageRemainder=\pagegoal
     \ifdim\PageRemainder=\maxdimen\PageRemainder=\vsize
     \else\advance\PageRemainder by -1\pagetotal\fi}

  \def\Rpt@{}\def\Rpt@@{}

  \long\def\HH#1\par{\par
  \SaveLastSkip\removelastskip\goodbreak
  \ifdim\LastSkip<30pt 
     \LastSkip 30pt
plus 3pt minus 2pt\fi
  \SetPageRemainder\advance\PageRemainder-\LastSkip
  \ifdim\PageRemainder<150pt
       \edef\Rpt@{remain = \the\PageRemainder\noexpand\\
                pagetotal=\the\pagetotal\noexpand\\
                           pagegoal=\the\pagegoal}%
          \fi
   \ifdim\PageRemainder<65pt 
       \ifdim\PageRemainder > 0pt
          \edef\Rpt@@{\noexpand\\
                      Had HH PageRemainder$<$\relax 65pt\noexpand\\
                      Hence forced break!}%
     \vskip 0pt plus .2\PageRemainder\eject 
    \fi\fi
    \vskip\LastSkip\Hrule 
    \pretolerance=10000\rightskip=0pt plus 3em
    \hangafter1 \hangindent=2.2em%
    \noindent
    \HHfont \unskip \Ednote{\Rpt@\Rpt@@}%
            \def\Rpt@{}\def\Rpt@@{}%
            \ignorespaces
            #1\par\rightskip=0pt\pretolerance=\StdPretolerance%
    \NoindentAfter
\tenpoint\rm%
     \medskip \vskip\ProcSkip}

  \long\def\HHH#1\par{\par%
  \SaveLastSkip\removelastskip\goodbreak
  \ifdim\LastSkip<\ProcSkip%
     \LastSkip\ProcSkip\fi
  \SetPageRemainder\advance\PageRemainder-\LastSkip
  \ifdim\PageRemainder<150pt
       \edef\Rpt@{remain = \the\PageRemainder\noexpand\\
                pagetotal=\the\pagetotal\noexpand\\
                           pagegoal=\the\pagegoal}%
       \fi
   \ifdim\PageRemainder<48pt  
        \ifdim\PageRemainder > 0pt
             \edef\Rpt@@{\noexpand\\
                      Had HHH PageRemainder$<$\relax48pt\noexpand\\
                      Hence forced break!}%
       \vskip 0pt plus .2\PageRemainder\eject 
      \fi\fi
   \vskip\LastSkip\par\noindent
   \HHHfont \unskip\Ednote{\Rpt@\Rpt@@}%
  \def\Rpt@{}\def\Rpt@@{}%
  \ignorespaces
   #1\unskip.\quad\rm\ignorespaces
   \ignorepars}

  \long\def\ignorepars#1\par{\def\Test{#1}%
     \ifx\Test\Empty\def\This{\ignorepars}%
        \else\def\This{\Test\par}\fi
           \This}
  \def\Empty{}

 \def\Abstract#1\par{\bgroup\Smallfonts\narrower\HHH #1\par}
 \def\endAbstract{\par\egroup}


 \def\ProcBreak{\par%
    \ifdim\lastskip<8pt%
    \removelastskip%
    \penalty-200\vskip\ProcSkip\fi}

 \def\th#1\par{\ProcBreak \noindent
   {\thfont\ignorespaces
    #1\unskip.}\it\itSpacing\kern.4em\ignorepars}

 \def\endth{\ProcBreak\rm\itSpacingOff }


 \def\pf#1\par{\ProcBreak %
    \noindent\pffont#1\unskip.\rm\itSpacingOff{\kern .7em}\ignorepars}

 \def\endpf{\medskip \ProcBreak } 

  \def\qedbox{\hbox{\vbox{
    \hrule width0.2cm height0.2pt
    \hbox to 0.2cm{\vrule height 0.2cm width 0.2pt
             \hfil\vrule height0.2cm width 0.2pt}
    \hrule width0.2cm height 0.2pt}\kern1pt}}

  \def\qed{\ifmmode\qedbox
    \else\unskip\ \hglue0mm\hfill\qedbox\ProcBreak\fi}

  \def \rk #1\par{\ProcBreak
     \noindent{\rkfont\ignorespaces #1\unskip.}%
     \rm\kern.6em\ignorepars}

  \def \endrk {\medskip\ProcBreak }

  \def \df #1\par{\ProcBreak
     \noindent{\dffont\unskip\ignorespaces #1\unskip.}%
     \rm\kern.6em\ignorepars}

  \def \enddf {\medskip\ProcBreak }

  \def \eg #1\par{\ProcBreak
     \noindent\egfont\unskip\ignorespaces #1\unskip.
     \rm\kern.6em\ignorepars}

  \newdimen\Overhang

   \def\MaxTag@#1#2#3#4#5{\setbox0=\hbox{#4\ignorespaces#2\unskip}%
     \dimen0=\wd0\advance\dimen0 by#3
     \ifdim\dimen0<#5\relax\dimen0=#5\fi
     \expandafter\edef\csname #1Hang\endcsname{\the\dimen0}}

 \def\MaxItemTag#1{\MaxTag@{Item}{#1}{.4em}{\ItemStyle}{\parindent}}%
 \def\MaxItemItemTag#1{%
        \MaxTag@{ItemItem}{#1}{.4em}{\ItemItemStyle}{\parindent}}
 \def\MaxNrTag#1{\MaxTag@{Nr}{#1}{.5em}{\NrStyle}{\parindent}}
 \def\MaxReferenceTag#1{%
        \MaxTag@{Reference}{[#1]}{.6em}{\ninerm}{\parindent}}
 \def\MaxFootTag#1{\MaxTag@{Foot}{#1}{.4em}{\ninerm}{\z@}}

  \def\SetOverhang@{\Overhang=.8\dimen0%
     \advance\Overhang by \wd0\relax
     \ifdim\Overhang>\hangindent\relax
       \advance\Overhang by .25\dimen0%
       \Ednote{Tag is pushing text.}\osumess{Tag is pushing text.}%
     \else\Overhang=\hangindent
     \fi}

   \def\Item#1{\par\noindent
      \hangafter1\hangindent=\ItemHang
      \setbox0=\hbox{\ItemStyle\ignorespaces#1\unskip}%
      \dimen0=.4em\SetOverhang@
      \rlap{\box0}\kern\Overhang\ignorespaces}

   \def\ItemItem#1{\par\noindent
      \hangafter1\hangindent=\ItemItemHang
      \setbox0=\hbox{\ItemItemStyle\ignorespaces#1\unskip}%
      \dimen0=.4em\SetOverhang@
      \advance\hangindent by \ItemHang
      \kern\ItemHang\rlap{\box0}%
      \kern\Overhang\ignorespaces}

  \def\Nr#1{\par\noindent\hangindent=\NrHang 
    \setbox0=\hbox{\NrStyle\ignorespaces#1\unskip}%
    \dimen0=.5em\SetOverhang@
    \rlap{\box0}\kern\Overhang
    \hangindent=\z@\ignorespaces}

   \newskip\Rosterskip\Rosterskip 1pt plus1pt 
   \def\Roster{\par\ifdim\lastskip<\Rosterskip\removelastskip\vskip\Rosterskip\fi
    \bgroup}
   \def\endRoster{\par\global\edef\LastSkip@{\the\lastskip}\removelastskip
       \egroup\penalty-50\LastSkip\LastSkip@\relax
       \ifdim\LastSkip<\Rosterskip\LastSkip\Rosterskip\fi
       \vskip\LastSkip}




 \def\cite#1{
    \def\nextiii@##1,##2\end@{{\frenchspacing\rm 
      \lBr\ignorespaces##1\unskip{\rm,~\ignorespaces##2}\rBr}}%
    \IN@0,@#1@%
    \ifIN@\def\next{\nextiii@#1\end@}\else
    \def\next{{\rm\lBr#1\rBr}}\fi\next}


   \def \Bib#1\par{%
       \par\removelastskip\SetPageRemainder
       \ifdim\PageRemainder < 97pt
        \ifdim\PageRemainder > 0pt
        \vfill\eject
       \fi\fi
    \ProcBreak \par\begingroup\parskip=0 pt%
    \goodbreak \vskip 15 pt plus 10 pt
    \noindent\null\hfill\Bibfont
      \ignorespaces #1\unskip\hfill\null\par 
    \frenchspacing \Smallfonts\rm
    \parskip=2.5 pt plus 1 pt minus.5pt%
    \nobreak\vskip 12pt plus 2pt minus2pt\nobreak
    \leftskip=0 pt \baselineskip=10.5pt}

 \def\ReferenceTagSlide{0em}
  \def\ReferenceTagGap{.5em}

  \def \rf#1{\par\noindent
     \hangafter1\hangindent=\ReferenceHang      
     \setbox0=\hbox{\ninerm[\ignorespaces#1\unskip]}%
     \dimen0=\ReferenceTagGap\SetOverhang@
     \rlap{\kern\ReferenceTagSlide\box0}%
     \kern\Overhang\ignorespaces}

  \def\ref#1\par#2\par#3\par#4\par{%
     \rf{#1}#2\unskip,\ #3\unskip,\
     #4\unskip.}

  \def\endBib{\par\endgroup\vskip 12pt minus 6pt }


  \long\def\Coordinates#1\endCoordinates{
 {\par\vskip4pt\def\\{\unskip, }\Coordfont\baselineskip10.5pt\noindent#1}}

 \def\pagecontents{
  \gdef\Pagetot@l{\pagetotal}
  \ifvoid\TRMargIns\else
    \rlap{\kern\hsize\kern10pt\vbox to 0pt{%
         \box\TRMargIns\vss}}\fi
  \ifvoid\topins\else\unvbox\topins\fi
   \dimen@=\dp\@cclv \unvbox\@cclv 
   \ifvoid\footins\else 
     \vskip\skip\footins
     \footnoterule
     \unvbox\footins\fi
   \ifr@ggedbottom \kern-\dimen@ \vfil \fi}


 \newcount\Ht 

 \def \Acc{\expandafter } 

 \def\swthat{\raise -1.1 ex\hbox{\sam$\widehat{}$}}
 \def\swttilde{\raise -1.2 ex\hbox{\sam$\widetilde{}$}}
 \def \overdot{{\raise .2 ex \hbox to 0pt {\hss\bf\smash{.}\hss}}}
 \def \overcircle{{\raise .1 ex \hbox to 0pt
    {\sam$\eightpoint\scriptstyle\hss\circ\hss$}}}

 \def \Mathaccent#1#2{{\sam 
  \setbox4=\hbox{$\vphantom{#2}$}
  \Ht=\ht4 
  \setbox5=\hbox{${#1}$}
  \setbox6=\hbox{${#2}$}
  \setbox7=\hbox to .5\wd6{}
  \copy7\kern .1\Ht \raise\Ht sp\hbox{\copy5}\kern-.1\Ht
  \copy7\llap{\box6}
  }}

  \def\SwtCheck #1{
        \ifmmode \check{#1}%
                \else \v {#1}%
                \fi}

 \def\barpartial {%
   \kern .17 em
    \overline {\kern -.17 em\partial\kern-.03 em}%
    \kern .03 em}

 
  \def\Overline#1{\setbox1=\hbox{\sam ${#1}$}%
      \ifdim \wd1 > 6pt
    \kern .11 em
    \overline {\kern -.11 em#1\kern-.14 em}
    \kern .14 em
  \else
    \kern .03 em
    \overline {\kern -.03 em#1\kern-.04 em}
    \kern .04 em
  \fi}

 \def\SOverline#1{\setbox1=\hbox{\sam ${#1}$}%
      \ifdim \wd1 > 7pt
    \kern .22 em
    \overline {\kern -.22 em#1\kern-.09 em}%
    \kern .09 em
  \else
    \kern .10 em
    \overline {\kern -.10 em#1\kern-.04 em}%
    \kern .04 em
  \fi}


 \def\Underline#1{\setbox1=\hbox{\sam ${#1}$}%
      \ifdim \wd1 > 6pt
    \kern .11 em
    \underline {\kern -.11 em#1\kern-.14 em}
    \kern .14 em
  \else
    \kern .03 em
    \underline {\kern -.03 em#1\kern-.04 em}
    \kern .04 em
  \fi}

 \def\SUnderline#1{\setbox1=\hbox{\sam ${#1}$}%
      \ifdim \wd1 > 7pt
    \kern .04 em
    \underline {\kern -.04 em#1\kern-.2 em}%
    \kern .2 em
  \else
    \kern .0 em
    \underline {\kern -.0 em#1\kern-.15 em}%
    \kern .15 em
  \fi}


 \def \Blackbox
   {\leavevmode\hskip .3pt \vbox
   {\hrule height 5pt\hbox{\hskip 4.5pt}}\hskip .5pt}

 \def \XX{\Blackbox\kern.5pt\Blackbox} 

  \def\.{.\kern1pt}

    \def\Hyphen{\edef\this{\the\hyphenchar\font}%
          \hyphenchar\font=-1\char\this\hyphenchar\font=\this}

 \ifx\undefined\text
  \def\text#1{\hbox{\rm #1}}\fi 



   \everymath{}  

  \def\PassMath@@{\aftergroup\AfterMath@} 

 \let\PassMath@\PassMath@@

 \def\AfterMath@{\futurelet\next\AfterMathMole@}

 \def\AfterMathMole@{
      \ifcat\next\space
          \def\this{}
      \else
      \ifcat\next\egroup %
        \def\this{\osumess{Handset mathsurround?? ---(see dollar brace)}}%
      \else
      \def\this{\AAfterMath@}
      \fi\fi
      \this}

 \def\hyphen@{-}
 \def\paren@{)}
 \def\apostr@{'}

 \def\MSC#1{\kern-.8\mathsurround#1\kern.8\mathsurround}

 \def\AAfterMath@#1{\def\Next{#1}
    \IN@0\Next @,.;:!?\relax @%
    \ifIN@\def\this{\MSC{\Next}}%
    \else
    \ifx\Next\hyphen@\def\this{\futurelet\next\AfterHyphen@}%
    \else
    \ifx\Next\paren@\def\this{#1}%
    \else 
    \ifx\Next\apostr@\def\this{#1}%
    \else \def\this{\osumess{Handset mathsurround??}%
                 #1}\fi\fi\fi\fi
    \this}

 \def\AfterHyphen@#1{\def\Next{#1}%
   \ifx\Next\hyphen@\def\this{--}\else
   \ifcat\next\space%
   \def\this{\kern-\mathsurround\kern.05em- \Next}\else
   \def\this{\kern-\mathsurround\kern.05em\Hyphen\Next}\fi\fi\this}

 \def\sam{\mathsurround=\z@\let\PassMath@\relax}  %
 \def\mas{\mathsurround=\StdMathsurround\let\PassMath@\PassMath@@}
 
 \def\Mas{\mathsurround=\StdMathsurround
                \everymath{\PassMath@}\let\PassMath@\PassMath@@}

 \def\m@th{\mathsurround=\z@\everymath{}}

 \def\m@@th{\mathsurround=\z@\everymath={}\let\m@th\relax}

\def\underbar#1{$\setbox\z@\hbox{#1}\dp\z@\z@
      \m@th \underline{\box\z@}$\relax}

\def\mathhexbox#1#2#3{\leavevmode
  \hbox{\m@@th$\m@th \mathchar"#1#2#3$}}

\def\dots{\relax\ifmmode\ldots\else$\m@th\ldots\,$\relax\fi}

\def\dotfill{\cleaders\hbox{\m@@th$\m@th \mkern1.5mu.\mkern1.5mu$}\hfill}
\def\rightarrowfill{$\m@th\mathord-\mkern-6mu%
  \cleaders\hbox{\m@@th$\mkern-2mu\mathord-\mkern-2mu$}\hfill
  \mkern-6mu\mathord\rightarrow$\relax}
\def\leftarrowfill{$\m@th\mathord\leftarrow\mkern-6mu%
  \cleaders\hbox{\m@@th$\mkern-2mu\mathord-\mkern-2mu$}\hfill
  \mkern-6mu\mathord-$\relax}

\def\downbracefill{$\m@th\braceld\leaders\vrule\hfill\braceru
  \bracelu\leaders\vrule\hfill\bracerd$\relax}
\def\upbracefill{$\m@th\bracelu\leaders\vrule\hfill\bracerd
  \braceld\leaders\vrule\hfill\braceru$\relax}

\def\angle{{\vbox{\m@@th\ialign{$\m@th\scriptstyle##$\crcr
      \not\mathrel{\mkern14mu}\crcr
      \noalign{\nointerlineskip}
      \mkern2.5mu\leaders\hrule height.34pt\hfill\mkern2.5mu\crcr}}}}

\def\big#1{{\m@@th\hbox{$\left#1\vbox to8.5\p@{}\right.\n@space$}}}
\def\Big#1{{\m@@th\hbox{$\left#1\vbox to11.5\p@{}\right.\n@space$}}}
\def\bigg#1{{\m@@th\hbox{$\left#1\vbox to14.5\p@{}\right.\n@space$}}}
\def\Bigg#1{{\m@@th\hbox{$\left#1\vbox to17.5\p@{}\right.\n@space$}}}
\def\n@space{\nulldelimiterspace\z@ \m@th}

\def\root#1\of{\setbox\rootbox\hbox{\m@@th$\m@th\scriptscriptstyle{#1}$}
  \mathpalette\r@@t}
\def\r@@t#1#2{\setbox\z@\hbox{\m@@th$\m@th#1\sqrt{#2}$\relax}
  \dimen@\ht\z@ \advance\dimen@-\dp\z@
  \mkern5mu\raise.6\dimen@\copy\rootbox \mkern-10mu \box\z@}

\def\mathph@nt#1#2{\setbox\z@\hbox{\m@@th$\m@th#1{#2}$}\finph@nt}

\def\mathsm@sh#1#2{\setbox\z@\hbox{\m@@th$\m@th#1{#2}$}\finsm@sh}

\def\@vereq#1#2{\lower.5\p@\vbox{\m@@th\baselineskip\z@skip\lineskip-.5\p@
    \ialign{$\m@th#1\hfil##\hfil$\crcr#2\crcr=\crcr}}}

\def\mathpalette#1#2{\sam\mathchoice{#1\displaystyle{#2}}%
  {#1\textstyle{#2}}{#1\scriptstyle{#2}}{#1\scriptscriptstyle{#2}}\mas}

\def\widehat#1{\setbox\z@\hbox{\sam$#1$}%
 \ifdim\wd\z@>\tw@ em\mathaccent"0\msbfam@5B{#1}%
 \else\mathaccent"0362{#1}\fi}
\def\widetilde#1{\setbox\z@\hbox{\sam$#1$}%
 \ifdim\wd\z@>\tw@ em\mathaccent"0\msbfam@5D{#1}%
 \else\mathaccent"0365{#1}\fi}

 \def\dots{\relax{}
  \ifmmode\def\thedots{\mdots@}\else\def\thedots{\tdots@}\fi %
  \thedots}

 \let\@ldeqno\eqno\let\@ldleqno\leqno
 \def\eqno{\everymath{}\@ldeqno} \def\leqno{\everymath{}\@ldleqno}

  \let\@ldeqalignno\eqalignno
  \def\eqalignno#1{\sam\@ldeqalignno{#1}\mas}
  \let\@ldeqalign\eqalign
  \def\eqalign#1{\sam\@ldeqalign{#1}\mas}

 \def\overrightarrow#1{\vbox{\m@th\ialign{##\crcr
      \rightarrowfill\crcr\noalign{\kern-\p@\nointerlineskip}
      $\hfil\displaystyle{#1}\hfil$\crcr}}}
 \def\overleftarrow#1{\vbox{\m@th\ialign{##\crcr
      \leftarrowfill\crcr\noalign{\kern-\p@\nointerlineskip}
      $\hfil\displaystyle{#1}\hfil$\crcr}}}
 \def\overbrace#1{\mathop{\vbox{\m@th\ialign{##\crcr\noalign{\kern3\p@}
      \downbracefill\crcr\noalign{\kern3\p@\nointerlineskip}
      $\hfil\displaystyle{#1}\hfil$\crcr}}}\limits}
 \def\underbrace#1{\mathop{\vtop{\m@th\ialign{##\crcr
      $\hfil\displaystyle{#1}\hfil$\crcr\noalign{\kern3\p@\nointerlineskip}
      \upbracefill\crcr\noalign{\kern3\p@}}}}\limits}

  \let\@ldmatrix\matrix
  \let\end@ldmatrix\endmatrix
  \def\matrix{\sam\@ldmatrix}
  \def\endmatrix{\end@ldmatrix\mas}
  \let\@ldgather\gather
  \let\end@ldgather\endgather
  \def\gather{\sam\@ldgather}
  \def\endgather{\end@ldgather\mas}
  \let\@ldalign\align
  \let\end@ldalign\endalign
  \def\align{\sam\@ldalign}
  \def\endalign{\end@ldalign\mas}
  \let\@ldaligned\aligned
  \let\end@ldaligned\endaligned
  \def\aligned{\sam\@ldaligned}
  \def\endaligned{\end@ldaligned\mas}
  \let\@ldtag\tag
  \def\tag{\sam\@ldtag}
   %

   \let\MinCDArrowWidth\minCDaw@




\newskip\insertskipamount\newskip\inserthardskipamount
\insertskipamount 6pt plus2pt 
\inserthardskipamount 6pt
\def\insertskip{\vskip\insertskipamount}
\newcount\SplitTest
\def\SetSplitTest{\SplitTest\insertpenalties
  \insert\topins{\floatingpenalty1}%
  \advance\SplitTest-\insertpenalties}
\def\midinsert{\par
 \SaveLastSkip\penalty-150\SetSplitTest\RestoreLastSkip
 \ifnum\SplitTest=-1
  \@midfalse\p@gefalse\else\@midtrue\fi\@ins}
\def\@ins{\par\begingroup\setbox\z@\vbox\bgroup%
  \vglue\inserthardskipamount}
\def\endinsert{\egroup 
  \if@mid \dimen@\ht\z@ \advance\dimen@\dp\z@
    \advance\dimen@\insertskipamount
    \advance\dimen@\pagetotal\advance\dimen@-\pageshrink
    \ifdim\dimen@>\pagegoal\@midfalse\p@gefalse\fi\fi
  \if@mid%
    \ifdim\lastskip<\insertskipamount\removelastskip\insertskip\fi
    \nointerlineskip\box\z@\penalty-200\insertskip
  \else%
    \SaveLastSkip
    \insert\topins{\penalty100 
    \splittopskip\z@skip
    \splitmaxdepth\maxdimen \floatingpenalty\z@
    \ifp@ge \dimen@\dp\z@
    \vbox to\vsize{\unvbox\z@\kern-\dimen@}
    \else \box\z@\nobreak\insertskip\fi}
    \RestoreLastSkip
   \fi\endgroup}


  \newcount\notenumber
  
  \def\note{\advance\notenumber by 1
    \footnote{\the\notenumber)}}

  \newbox\footbox

  \def\footnote#1{\let\@sf\empty
    \ifhmode\edef\@sf{\spacefactor\the\spacefactor}\/\fi
    \sam${}^{\fam0 #1}$\@sf\vfootnote{#1}}%

  \def\vfootnote#1{\insert\footins\bgroup
     \interlinepenalty100 \splittopskip=1pt
     \floatingpenalty=20000
     \leftskip=0pt\rightskip=0pt%
     \parindent=.3em
     \Smallfonts\rm
     \FootItem@{#1}
     \futurelet\next\fo@t}

  \def\FootItem@#1{\par\hangafter1\hangindent=\FootHang
     \setbox0=\hbox{\ignorespaces#1\unskip}%
     \dimen0=.4em\SetOverhang@
     \noindent\rlap{\box0}\kern\Overhang\ignorespaces}


  \def\fo@t{\ifcat\bgroup\noexpand\next \let\next\f@@t
    \else\let\next\f@t\fi \next}
  \def\f@@t{\bgroup\aftergroup\@foot\let\next}
  \def\f@t#1{\baselineskip=10pt\lineskip=1pt
            \lineskiplimit=0pt #1\@foot}%
  \def\@foot{
        \hbox{\vrule height0pt depth5pt width0pt}
        \egroup}
  \skip\footins=12 pt plus 0pt minus 0pt 
  \count\footins=1000 
  \dimen\footins=8in 



 \def\osumess#1{\EdSpider{\immediate\write16{Line \the\inputlineno: #1}}}%
 \def\HideEdStuff{\gdef\EdSpider##1{}}

 \font\BigSym=cmmi10 scaled \magstep 4

 \def\change{\InLMargin{\hbox{\BigSym \char63\kern10pt}}}

 \def\beginchange{\InLMargin{\hbox{\sam\twelvepoint$\heartsuit$\kern10pt}}}

 \def\endchange{\InLMargin{\hbox{\sam\twelvepoint$\spadesuit$\kern10pt}}}

 \def\InLMargin#1{\strut\vadjust{%
     \kern-\strutdepth
     \vtop to \strutdepth{%
         \baselineskip\strutdepth
         \llap{\sam$\smash{\hbox{\EdSpider{#1}}}$}\null}}}

 \def\strutdepth{\dp\strutbox}
 \def\strutheight{\ht\strutbox}

 \def\NoteInRMargin#1{\strut\vadjust{%
     \kern-1.001\strutdepth
     \vtop to \strutdepth{%
       \baselineskip\strutdepth
       \vss\rlap{\ninepoint\unskip\hskip\hsize
         \vtop to 0pt{%
           \hsize=16em\hfuzz=\hsize
           \leftskip=10pt%
           \rightskip=0pt plus 10000pt%
           \baselineskip=9.8pt\lineskip=.2pt%
           \let\\\break
           \noindent\EdSpider{#1}\vss}%
                \kern10pt}\hbox{}}
       }}

 \def\ednote#1{\NoteInRMargin{\tentt #1}}

 \def\cbar{\InLMargin{%
      \dimen0=\strutdepth\advance\dimen0 by \lineskip
      \vrule width 3pt
      height \strutheight depth \dimen0 \kern
      3pt}}

 \def\ccbar{\InLMargin{%
      \dimen0=2\strutdepth\advance\dimen0 by 2\lineskip
      \vrule width 3pt
        height 3\strutheight depth \dimen0 \kern
      3pt}}

 \newinsert\TRMargIns
 \dimen\TRMargIns=\maxdimen

  \def\Ednote#1{\insert\TRMargIns{%
       \vbox to 0pt{\hsize=140pt\hfuzz=\hsize
           \leftskip=6pt%
           \rightskip=0pt plus 10000pt%
           \baselineskip=9.8pt\lineskip=.2pt%
           \let\\\break
           \SetPageRemainder
           \vglue540pt\vglue-\PageRemainder
           \noindent\EdSpider{\tentt #1}\vss}%
       \smallskip}}

 \def\KillEdStuff{\def\ednote##1{}\def\Ednote##1{}%
      \let\change\relax\let\beginchange\relax\let\endchange\relax
       \let\cbar\relax\let\ccbar\relax}


  \topskip=12pt
  \newskip\StdBaselineskip 
  \StdBaselineskip 12pt
  \lineskip=1.1pt
  \lineskiplimit=.8pt
  \widowpenalty=10000 
  \clubpenalty=10000  
  \abovedisplayskip=6pt plus 1pt minus 1pt
  \abovedisplayshortskip=3pt plus 1.5pt
  \belowdisplayskip=6pt plus 1pt minus 1pt
  \belowdisplayshortskip=5pt plus 1pt minus 1pt
  \hfuzz=1.5pt   

  \def\StdPretolerance{100}
  \tolerance=\StdPretolerance

  \newdimen\StdMathsurround
  \StdMathsurround=1.5pt 
  \mathsurround=\StdMathsurround
  \Mas                   

   \def\prose{\relax\hbox{\kern.6\StdMathsurround}}
  
  \def\StdParskip{0pt}    
  \parskip=\StdParskip
  \parindent=0.5cm
 

  \def\Times{ptmr  } 
  \def\TimesI{ptmri  } 
  \def\TimesB{ptmb  }
  \def\TimesBI{ptmbi  }
  \def\HelveticaN{phvrrn }

  =\Times at 10bp
  =\TimesB at 10bp
  \font\tenit=\TimesI at 10bp
  =\TimesBI at 10bp

  \font\tenmrm=cmr10  


    =\Times at 9bp 
    \font\nineit=\TimesI at 9bp 
    =\TimesB at 9bp 
    =\TimesBI at 9bp 

    =\HelveticaN at 9bp 


  =\Times at 12bp
  \font\twelveit=\TimesI at 12bp
  =\TimesB at 12bp


  \font\titleit=\TimesI at 14.4bp
  =\TimesB at 14.4bp

 \SetAuthorHead{AuthorHead} 
 \SetTitleHead{TitleHead}  


  \def\lBr{\raise.125ex\hbox{[\kern.1125ex}}
  \def\rBr{\raise.125ex\hbox{\kern.1125ex]}}

 \setbox\footbox=\hbox{\Smallfonts 2)~}



  \bgroup
  \catcode`\@=11 
  \gdef\itSpacing{%
     \xspaceskip=.31em plus.1em minus.05em \sfcode `f=2001
     \itWarning@\let\itWarning@\itWarning@@}
  \gdef\itSpacingOff{%
     \xspaceskip=0pt \sfcode `f=1000
     \let\itWarning@\relax}
   \global\let\itWarning@\relax
  \gdef\itWarning@@{\errmessage{%
  Special italic spacing already in force
  (you have probably omitted an ``endth'').
  See itSpacing macro in osuPSfnt.sty
         }}
  \egroup

 \fontdimen1\titlebf=0.0pt
 \fontdimen2\titlebf=3.6135pt
 \fontdimen3\titlebf=2.8908pt
 \fontdimen4\titlebf=1.44539pt
 \fontdimen5\titlebf=6.64882pt
 \fontdimen6\titlebf=14.45398pt
 \fontdimen7\titlebf=1.60439pt

 \fontdimen1\tenbi=0.26794pt
 \fontdimen2\tenbi=2.50937pt
 \fontdimen3\tenbi=2.00749pt
 \fontdimen4\tenbi=1.00374pt
 \fontdimen5\tenbi=4.59717pt
 \fontdimen6\tenbi=10.03749pt
 \fontdimen7\tenbi=1.11415pt

 \fontdimen1\twelverm=0.0pt
 \fontdimen2\twelverm=3.01125pt
 \fontdimen3\twelverm=2.409pt
 \fontdimen4\twelverm=1.2045pt
 \fontdimen5\twelverm=5.39615pt
 \fontdimen6\twelverm=12.045pt
 \fontdimen7\twelverm=1.33699pt

 \fontdimen1\twelveit=0.27731pt
 \fontdimen2\twelveit=3.01125pt
 \fontdimen3\twelveit=2.409pt
 \fontdimen4\twelveit=1.2045pt
 \fontdimen5\twelveit=5.37207pt
 \fontdimen6\twelveit=12.045pt
 \fontdimen7\twelveit=1.33699pt

 \fontdimen1\twelvebf=0.0pt
 \fontdimen2\twelvebf=3.01125pt
 \fontdimen3\twelvebf=2.409pt
 \fontdimen4\twelvebf=1.2045pt
 \fontdimen5\twelvebf=5.5407pt
 \fontdimen6\twelvebf=12.045pt
 \fontdimen7\twelvebf=1.33699pt

 \fontdimen1\tenrm=0.0pt
 \fontdimen2\tenrm=2.50937pt
 \fontdimen3\tenrm=2.00749pt
 \fontdimen4\tenrm=1.00374pt
 \fontdimen5\tenrm=4.49678pt
 \fontdimen6\tenrm=10.03749pt
 \fontdimen7\tenrm=1.11415pt

 \fontdimen1\tenit=0.27731pt
 \fontdimen2\tenit=2.50937pt
 \fontdimen3\tenit=2.00749pt
 \fontdimen4\tenit=1.00374pt
 \fontdimen5\tenit=4.47672pt
 \fontdimen6\tenit=10.03749pt
 \fontdimen7\tenit=1.11415pt

 \fontdimen1\tenbf=0.0pt
 \fontdimen2\tenbf=2.50937pt
 \fontdimen3\tenbf=2.00749pt
 \fontdimen4\tenbf=1.00374pt
 \fontdimen5\tenbf=4.61723pt
 \fontdimen6\tenbf=10.03749pt
 \fontdimen7\tenbf=1.11415pt

 \fontdimen1\ninerm=0.0pt
 \fontdimen2\ninerm=2.25842pt
 \fontdimen3\ninerm=1.80673pt
 \fontdimen4\ninerm=0.90337pt
 \fontdimen5\ninerm=4.0471pt
 \fontdimen6\ninerm=9.03374pt
 \fontdimen7\ninerm=1.00273pt

 \fontdimen1\nineit=0.27731pt
 \fontdimen2\nineit=2.25842pt
 \fontdimen3\nineit=1.80673pt
 \fontdimen4\nineit=0.90337pt
 \fontdimen5\nineit=4.02904pt
 \fontdimen6\nineit=9.03374pt
 \fontdimen7\nineit=1.00273pt

 \fontdimen1\ninebf=0.0pt
 \fontdimen2\ninebf=2.25842pt
 \fontdimen3\ninebf=1.80673pt
 \fontdimen4\ninebf=0.90337pt
 \fontdimen5\ninebf=4.15552pt
 \fontdimen6\ninebf=9.03374pt
 \fontdimen7\ninebf=1.00273pt


 \newcount\MaxSpaceFactor
 \MaxSpaceFactor=3000 

 \def\ItemStyle{\rm}
 \def\NrStyle{\rm}
 \def\ItemItemStyle{\rm}

 \MaxItemTag{(iii)}
 \MaxItemItemTag{(iii)}
 \MaxNrTag{(2)}
 \MaxFootTag{2)}
 \def\ReferenceHang{30pt}

 \catcode`\@=\active


\loadbold

=\Times  
=\Times scaled750
=\Times scaled650
\font\rms=\Times scaled 920 

=\TimesBI scaled 860
=\TimesI scaled 860

\textfont0=\rrm  
\scriptfont0=\erm 
\scriptscriptfont0=\srm

\def\Augment#1#2{%
    \toks0\expandafter{#1}\toks2{#2}%
    \edef#1{\the\toks0\the\toks2}}

 \font\twelverma=\Times  scaled 1200
 \font\tenrma=\Times  scaled 1000
 \font\ninerma=\Times scaled 920
 =\Times scaled 840
 \font\sevenrma=\Times scaled 760
 =\Times scaled 680
 \font\fiverma=\Times scaled 600

 \Augment\tenpoint{%
  \textfont0=\tenrma  \scriptfont0=\sevenrma  
  \scriptscriptfont0=\fiverma  }

 \Augment\ninepoint{%
  \textfont0=\ninerma  \scriptfont0=\sevenrma 
  \scriptscriptfont0=\fiverma}

 \Augment\twelvepoint{%
  \textfont0=\twelverma  \scriptfont0=\ninerma  
  \scriptscriptfont0=\sevenrma}

\mathsurround=1pt
\hsize=13.45truecm
\vsize=19.5truecm
\hoffset=1.25truecm
\voffset=2truecm
\advance\baselineskip by 2pt

\predefine\til{\~}
\def\~#1{\relax\ifmmode\widetilde{#1}\else\til{#1}\fi}

\redefine \le{\leqslant}
\redefine \ge{\geqslant}
\define \wt#1{\mathaccent"0365{#1}}
\define \wh#1{\mathaccent"0362{#1}}

\define \iss{\,\Mathaccent{\raise -.8 ex\hbox{$\widetilde{}$\kern.1em}}\rightarrow\,}

\define \minn{\operatorname{\fam0 min\,}}

\define \tpp{\mathop{\fam0 top}}
\define \ab{\mathop{\fam0 ab}}

\define \kr{\mathop{\fam0 ker}}

\define \chr{\mathop{\fam0 char}\,}

\define \res{\operatorname{\fam0 res}}

\define \Tr{\operatorname{\fam0 Tr\,}}
\define \Tor{\operatorname{\fam0 Tors}}

\define \Gal{\mathop{\fam0 Gal}}

\Mas
\HideEdStuff
\rm 
 

\def\issn{{\nineit ISSN 1464-8997 (on line) 1464-8989 (printed)}}

\def\gtp{{\nineit Published 10 December 2000: \ \copyright\ Geometry \& 
Topology Publications}}

\def\gtv3{{\nineit Geometry \& Topology Monographs, Volume 3 (2000) --
Invitation to higher local fields}}


\def\lione
{{\rms Geometry \& Topology Monographs}}

\def \litwo{{\rms Volume 3: Invitation to higher local fields
}} 

\def\tinfo #1.#2.#3-#4
{{
\noindent  {\lione} \hfill 
\par 
\vskip-1.5pt
\noindent {\litwo} \hfill
\par 
\vskip-1,5pt
\noindent {\rms Part #1, section #2, pages #3--#4} \hfill
\vskip24pt 
}}

\def\tinfos #1.#2.#3-#4
{{
\noindent  {\lione} \hfill 
\par 
\vskip-1.5pt
\noindent {\litwo} \hfill
\par 
\vskip-1.5pt
\noindent {\rms Pages #3--#4} \hfill
\vskip24pt 
}}

\def\tinfoi #1
{{
\noindent  {\lione} \hfill 
\par 
\vskip-1.5pt
\noindent {\litwo} \hfill
\par 
\vskip-1.5pt
\noindent {\rms Pages iii--xi: Introduction and contents} \hfill
\vskip26pt 
}}


  \def\titlepagehead{\hfil}

  \newif\iftitlepage\titlepagefalse
  \newif\ifblankpage\blankpagefalse
  \def\makeheadline{
     \ifblankpage{}\else%
     \iftitlepage
\vbox{\line{\vbox to 8.5pt{}
\ninerm
\copy\HLinebox \hfill
\hglue5mm\ninebf\folio 
\titlepagehead}}%
      \else
\vbox{\ifodd\pageno\rightheadline\else\leftheadline\fi}%
      \fi\vskip 12pt\fi}%
     \def\rightheadline{\line{\vbox to 8.5pt{}%
      \ninerm
\copy\TitleBox \hfill
\hglue5mm\ninebf\folio}}%
     \def\leftheadline{\line{\vbox to 8.5pt{}%
        \unskip\ninerm\unskip\ninebf\folio\hglue5mm
 \hfill \copy\AuthorBox
}}

 \footline={\ifblankpage{}\else
\iftitlepage\ninepoint\sam\hfill
\line{\vbox to 8.5pt{}
\copy\TFLinebox
\hfill
\hglue5mm 
}
            \else
\ninepoint\sam\hfill
\line{\vbox to 8.5pt{}
\copy\FLinebox
\hfill 
\hglue5mm
}
\hfil\fi\global\titlepagefalse\fi}

\def\blankpage{{\blankpagetrue\noindent\hbox to 10pt{\hss}\vfill
\pagebreak}}

\tenpoint\rm 
 

\pageno=61

\tinfo I.6.61-74

\SetTFLinebox{\gtp }
\SetFLinebox{\gtv3 }
\SetHLinebox{\issn}

\H 6. Topological Milnor $K$-groups of higher local fields 

Ivan Fesenko

\SetAuthorHead{I. Fesenko}
\SetTitleHead{Part I. Section 6.  Topological Milnor $K$-groups 
of higher local fields \qquad\qquad}

Let $F=K_n,\dots, K_0=\Bbb F_q$ be an $n$-dimensional local field.
We use the notation of section~1.

In this section we describe properties of
certain quotients $K^{\tpp}(F)$ of the Milnor $K$-groups of $F$
by using in particular topological considerations.
This is an updated and simplified summary of relevant results in 
\cite{F1--F5}. 
Subsection~6.1 recalls well-known results on
$K$-groups of classical local fields.
Subsection~6.2 discusses so called sequential topologies
which are important for the description of 
subquotients of  $K^{\tpp}(F)$ in terms of a simpler objects
endowed with sequential topology 
(Theorem~1 in 6.6 and Theorem~1 in~7.2 of section~7). 
Subsection~6.3 introduces $K^{\tpp}(F)$,
6.4 presents very useful pairings (including Vostokov's symbol
which is discussed in more detail in section~8),
subsection~6.5--6.6 describe the structure of $K^{\tpp}(F)$
and 6.7 deals with the quotients $K(F)/l$;
finally, 6.8 presents various properties of the norm map on $K$-groups.
Note that subsections 6.6--6.8 are not required for understanding 
 Parshin's class field theory in section~7.

\HH 6.0. Introduction

Let $A$ be a commutative ring and let $X$ be an $A$-module endowed with some
topology. A set $\{ x_i\}_{i\in I}$ of elements of $X$
is called a set of {\it topological generators} of $X$
if the sequential closure of the submodule of $X$ generated by this set
coincides with $X$. A set of topological generators is called
a {\it topological basis} if for every $j\in I$ and every non-zero $a\in A$
$ax_j$ doesn't belong to the sequential closure of the submodule
generated by $\{x_i\}_{i\not=j}$.

Let $I$ be a countable set. If $\{x_i\}$ is set of topological generators
of $X$ then every element $x\in X$ can be expressed as a convergent sum
$\sum a_ix_i$ with some $a_i\in A$ (note that it is not necessarily the case that
for all $a_i\in A$ the sum $\sum a_ix_i$ converges).
This expression is unique if $\{x_i\}$ is a topological basis of $X$;
then provided  addition in $X$ is sequentially continuous, we get 
$\sum a_ix_i+\sum b_ix_i=\sum (a_i+b_i)x_i$. 

Recall that in the  one-dimensional case the group of principal units
$U_{1,F}$ is a multiplicative  $\Bbb Z_p$-module
with finitely many topological generators if $\chr(F)=0$ and infinitely many topological generators
if $\chr(F)=p$ (see for instance \cite{FV, Ch. I \S6}).
This representation of $U_{1,F}$ and
a certain specific choice of its generators 
is quite important if one wants to deduce 
the Shafarevich and Vostokov explicit formulas for the Hilbert symbol
(see section~8).

Similarly,  the group $V_F$ of principal units
of an $n$-dimensional local field $F$ is topologically generated
by $1+\theta t_n^{i_n}\dots t_1^{i_1}$, $\theta\in \mu_{q-1}$ (see subsection~1.4.2). 
This leads to a natural suggestion to 
endow the  Milnor $K$-groups of $F$ with an appropriate topology and 
use the  sequential convergence to simplify calculations in $K$-groups. 

On the other hand, 
the reciprocity map $$\Psi _{F}\colon K_{n}(F)\rightarrow  
\Gal(F^{\ab}/F)$$ is not injective
in general, in particular 
$ \kr  (\Psi _{F})\supset \bigcap_{l\ge 
1}lK_{n}(F)\not=0$.
So  the Milnor $K$-groups are too large from the point of view of class field theory, and one can  pass to the quotient
$K_n(F)/\bigcap_{l\ge 
1}lK_{n}(F)$ without loosing any arithmetical information on $F$.
The latter quotient coincides with $K_n^{\tpp}(F)$  
(see subsection~6.6) which is defined in subsection~6.3 as the quotient
of $K_n(F)$ by the intersection $\Lambda_n(F)$ of all neighbourhoods of 0
in $K_n(F)$ with respect to a certain topology. 
The existence theorem in class field theory uses the topology to characterize  norm  subgroups $N_{L/F}K_n(L)$ 
of finite abelian extensions $L$ of $F$ as open subgroups of finite index 
of $K_n(F)$ (see subsection~10.5).   
As a corollary of the existence theorem in 10.5
one obtains that in fact 
$$\bigcap_{l\ge 1} lK_{n}(F)=\Lambda_n(F)= \kr  (\Psi _{F}).$$  
However, the class of  open subgroups of finite index of $K_n(F)$  can be defined without introducing the topology 
on $K_n(F)$,   
see the paper of Kato in this volume which presents a different approach.

\HH 6.1. $K$-groups of one-dimensional local fields

The structure of the Milnor $K$-groups of 
a one-dimensional local field $F$
is completely known.

Recall that using  the Hilbert symbol
and  multiplicative $\Bbb Z_p$-basis of the group of principal units of $F$
one obtains that
$$K_2(F)\simeq \Tor K_2(F)\oplus mK_2(F), 
\qquad\text{ 
where $m=|\Tor F^*|$, $\Tor K_2(F)\simeq \Bbb Z/m$ } $$
and $mK_2(F)$ is an uncountable uniquely divisible group 
(Bass, Tate, Moore, Merkur'ev; see for instance \cite{FV, Ch. IX \S4}).
The groups $K_m(F)$ for $m\ge 3$ are uniquely divisible uncountable groups
(Kahn \cite{Kn}, Sivitsky \cite{FV, Ch. IX \S4}).

\HH 6.2. Sequential topology

We need  slightly different topologies from the topology of $F$ and $F^*$ 
introduced in section~1. 
 
\df Definition

 Let $X$ be a topological space with topology $\tau 
$.  Define its {\it sequential saturation} $\lambda $: 

\noindent a subset $U$ of $X$ is
open with respect to $\lambda $ if for every $\alpha \in U$ and a convergent
(with respect to  $\tau $) sequence $X\ni \alpha _{i}$ to $\alpha $ almost all $\alpha
_{i} $ belong to $U$.  Then $\alpha _{i} @>>\tau>\alpha \Leftrightarrow \alpha _{i}  @>>\lambda>\alpha $.
\enddf

Hence the sequential saturation is the strongest topology which has the same
convergent sequences  and their limits as the original one. 
For a very elementary introduction to sequential topologies
see  \cite{S}. 

\df Definition

For  an $n$-dimensional local field $F$  denote by  $\lambda $ the sequential 
saturation of the topology on $F$ defined in section 1.
\enddf

The topology $\lambda $ is different from the old topology on $F$ defined in
section 1 for $n \ge 2$: for example, 
if $
F= \Bbb F_{p}\left( \left( t_{1}\right) \right) \left( \left(
t_{2}\right) \right)$ then
$Y=F\setminus \bigl\{
t_{2}^{i}t_{1}^{-j}+t_{2}^{-i}t_{1}^{j} : i,j \ge 1\bigr\} $ 
is open with respect to  $%
\lambda $ and is not open with respect to  the topology of $F$ defined in section 1. 

Let $\lambda _{* }$ on $F^{* }$ be the sequential saturation of the topology $\tau$
on $F^{*}$ defined in section 1.
It is a shift invariant topology. 

If $n=1$, the restriction of $\lambda _{* }$ on $V_F$
coincides with the induced from $\lambda$.

The following properties of $\lambda$ ($\lambda_*$) are similar to those in section~1
and/or can be  proved by induction on dimension. 
\rk Properties 

\Roster

\Item{(1)} $\alpha _{i},\beta _{i} @>>\lambda >0\Rightarrow
\alpha _{i}-\beta _{i} @>>\lambda > 0$;

\Item{(2)} $\alpha _{i},\beta _{i} @>>\lambda_* >1\Rightarrow \alpha _{i}\beta _{i}^{-1}  
@>>\lambda_* >1$;

\Item{(3)} for every $\alpha _{i}\in U_{F}$, $\alpha _{i}^{p^{i}}
@>>\lambda_* >1$;

\Item{(4)}  multiplication is not continuous in general
with respect to $\lambda_*$;  

\Item{(5)} every fundamental sequence with respect to $\lambda$ 
(resp. $\lambda_*$) converges; 

\Item{(6)} $V_F$ and $F^{*\, m}$ are closed subgroups of $F^*$
for every $m\ge 1$; 

\Item{(7)} The intersection of all
open subgroups of finite index containing a closed subgroup $H$ coincides
with $H$.
\endRoster  
\endrk

\df Definition

For topological spaces $X_{1},\dots ,X_{j}$ define the $* $-product topology
on $X_{1}\times \dots \times X_{j}$ as the sequential saturation of the product
topology.
\enddf

\HH 6.3. $K^{\tpp}$-groups

\df Definition

Let $\lambda _{m}$ be the strongest
topology on $K_{m}(F)$ such that
 subtraction in $K_m(F)$ and the natural map 
$$\varphi \colon (F^{* })^{m}\rightarrow K_{m}(F),
\quad \varphi(\alpha_1,\dots,\alpha_m)=\{(\alpha_1,\dots,\alpha_m\}$$
are sequentially continuous.
Then the
topology $\lambda _{m}$ coincides with its sequential saturation.
Put $$\Lambda _{m}(F)=\bigcap \text{ open neighbourhoods of 0}.$$
It is straightforward to see that $\Lambda_m(F)$ is a  subgroup of $K_m(F)$.
\enddf

\rk Properties

\Roster 

\Item{(1)} $\Lambda _{m}(F)$ is 
closed: indeed $\Lambda _{m}(F)\ni x_{i}\rightarrow x$
implies that  $x=x_{i}+y_{i}$ with 
$y_{i}\rightarrow 0$, so $x_{i},y_{i}\rightarrow 0$,  
hence $%
x=x_{i}+y_{i}\rightarrow 0$, so $x\in \Lambda _{m}(F)$. 

\Item{(2)} Put $VK_{m}(F)=\left\langle \left\{ V_F\right\} \cdot K_{m-1}(F)\right\rangle $ ($V_F$ is defined in subsection 1.1).
Since the topology with $VK_m(F)$ and its shifts as a system of fundamental neighbourhoods  satisfies 
two conditions of the previous definition, one obtains  that 
 $\Lambda 
_{m}(F)\subset VK_{m}(F)$. 

\Item{(3)} $\lambda _{1}=\lambda_*$.
\endRoster 

Following the original approach of Parshin \cite{P1}  
introduce now the following:

\df Definition

Set $$K_{m}^{  {\tpp}}(F)=K_{m}(F)/\Lambda _{m}(F)$$ 
and endow it with the quotient
topology of $\lambda_m$ which we denote by the same notation. 
\enddf

This new group $K_{m}^{  {\tpp}}(F)$ is sometimes called
the {\it topological Milnor $K$-group} of $F$.

If $\chr(K_{n-1})=p$ then $K_1^{\tpp} =K_{1}$.  

If $\chr
(K_{n-1})=0$ then $K_{1}^ {\tpp}(K)\neq K_{1}(K)$, since $1+ \Cal M%
_{K_{n}}$ (which is uniquely divisible)  is a subgroup of $\Lambda_1(K)$.

\HH 6.4. Explicit pairings

Explicit pairings of  the Milnor $K$-groups of $F$ 
are quite useful if one wants to study the structure
of $K^{\tpp}$-groups.

The general method is as follows.
Assume that there is a  pairing
$$\langle\,\,,\,\,\rangle\colon A\times B\to \Bbb Z/m$$
of two $\Bbb Z/m$-modules $A$ and $B$.
Assume that $A$ is endowed with a topology
with respect to which it has topological generators
$\alpha_i$ where $i$ runs over elements of a totally ordered countable set $I$.
Assume that
for every $j\in I$ there is an element
$\beta_j\in B$ such that
$$\langle \alpha_j,\beta_j\rangle =1 \mod m,\qquad
\langle \alpha_i,\beta_j\rangle =0 \mod m \quad
\text{for all $i>j$}.
$$
Then if a convergent sum $\sum c_i\alpha_i$ is equal to $0$, 
assume that there is a minimal $j$ with non-zero $c_j$
and deduce that
$$0=\sum c_i \langle\alpha_i ,\beta_j\rangle =c_j,$$
a contradiction.
Thus, $\{\alpha_i\}$ form a topological basis of $A$.

If, in addition, for every $\beta\in B\setminus\{0\}$
there is an $\alpha\in A$ such that
$\langle\alpha,\beta\rangle\not=0$, then the pairing $\langle\,\,,\,\,\rangle$ is obviously non-degenerate.

Pairings listed below satisfy the assumptions above
and therefore can be applied  
to study the structure of quotients of  the Milnor $K$-groups of $F$.

\HHH 6.4.1. ``Valuation map''

\phantom{}\par

Let $\partial\colon K_r(K_s)\to K_{r-1}(K_{s-1})$ be the border homomorphism (see for example \cite{FV, Ch. IX \S 2}). Put 
$${\goth v}={\goth v}_F\colon K_{n}(F)
@>\partial>> K_{n-1}(K_{n-1})@>\partial>>\dots @>\partial>>
K_0(K_0)=\Bbb Z, \quad {\goth v}(\{t_1,\dots,t_n\})=1$$
for a system of local parameters $t_1,\dots, t_n$ of $F$.
The valuation map $\goth v$ doesn't depend on the choice of a system of local parameters.

\HHH 6.4.2. Tame symbol

\phantom{}\par

Define 
$$t\colon K_{n}(F)/(q-1)\times F^{* }/F^{*\, q-1 }@>>>
K_{n+1}(F)/(q-1)@>>> \Bbb F_{q}^{\,* }\to \mu_{q-1},\quad q=|K_{0}|$$ by 
$$
K_{n+1}(F)/(q-1)
@>\partial>>K_{n}(K_{n-1})/(q-1)@>\partial>>\dots
@>\partial>>K_1(K_0)/(q-1)=
 \Bbb F_{q}^{\,* }\to \mu_{q-1}.$$
Here the map $\Bbb F_{q}^{\,* }\to \mu_{q-1}$ is given by
taking multiplicative representatives.

An explicit formula for this symbol (originally asked for in \cite{P2} and suggested in \cite{F1}) is simple:
let $t_1,\dots,t_n$ be a system of local parameters of $F$ and let $\bold{v}
=(v_1,\dots,v_n)$
be the associated valuation of rank $n$ (see section 1 of this volume). 
For elements $\alpha_1,\dots,\alpha_{n+1}$ of $F^*$  the value 
$t(\alpha_1,\alpha_2,\dots,\alpha_{n+1})$ is equal to 
the $(q-1)$th root of unity whose residue is equal to the residue
of 
$$\alpha_1^{b_1}\dots\alpha_{n+1}^{b_{n+1}}(-1)^b$$ 
in the last residue field $\Bbb F_q$, 
where $b=\sum_{s,i<j}v_s(b_i)v_s(b_j)b_{i,j}^s$, 
$b_j$ is the determinant
of the matrix obtained by cutting off  
the $j$th column of the
matrix $A=(v_i(\alpha_j))$  with the sign $(-1)^{j-1}$, and 
$b_{i,j}^s$ is the determinant of the matrix obtained by cutting off   
the $i$th and $j$th columns and $s$th row of $A$.

\HHH 6.4.3. Artin--Schreier--Witt pairing in characteristic $p$

\phantom{}\par

Define, following \cite{P2}, the pairing  
$$
(\;,\;]_{r}\colon K_{n}(F)/p^{r}\times W_{r}(F)/( {\bold F}-1)W_{r}(F)\rightarrow
W_{r}( \Bbb F_{p})\simeq  \Bbb Z/p^{r}
$$
by (${\bold F}$ is the map defined in the section, Some Conventions) 
$$
(\alpha _{0},\dots ,\alpha _{n},(\beta _{0},\dots ,\beta _{r})]_{r}=\Tr%
_{K_{0}/ \Bbb F_{p}}\quad (\gamma _{0},\dots ,\gamma _{r})
$$
where the $i$th ghost component $\gamma ^{(i)}$
is given by $\res_{K_{0}}$ $(\beta
^{(i)}\alpha _{1}^{-1}d\alpha _{1}\wedge \dots \wedge \alpha
_{n}^{-1}d\alpha _{n})$.

For its properties see \cite{P2, sect. 3}.
In particular, 
\Roster
\Item{(1)} for $x\in K_{n}(F)$ 
$$(x,\bold V(\beta_0,\dots,\beta_{r-1})]_r=\bold V(x,(\beta_0,\dots,\beta_{r-1})]_{r-1}$$
where as usual for a field $K$ 
$$\bold V\colon W_{r-1}(K)\to W_r(K), \quad \bold V(\beta_0,\dots,\beta_{r-1})=(0,\beta_0,\dots,\beta_{r-1});$$

\Item{(2)} for $x\in K_{n}(F)$ 
$$(x,\bold A(\beta_0,\dots,\beta_{r})]_{r-1}=\bold A(x,(\beta_0,\dots,\beta_{r})]_{r}$$
where for a field $K$ 
$$\bold A\colon W_{r}(K)\to W_{r-1}(K), \quad 
\bold A(\beta_0,\dots,\beta_{r-1},\beta_r)=(\beta_0,\dots,\beta_{r-1}).$$ 

\Item{(3)}If $\Tr\theta _{0}=1$ then  
$ 
\bigl(\{ t_{1},\dots ,t_{n}\} ,\theta _{0}\bigr]_{1}=1$.
If $i_l$ is prime to $p$ then 
$$\bigl(\{ 1+\theta t_{n}^{i_{n}}\dots t_{1}^{i_{1}},t_1,\dots  ,\widehat{t_l},\dots, t_n\} ,\theta _{0}\theta ^{-1}i_{l}^{-1}t_{1}^{-i_{1}}\dots
t_{n}^{-i_{n}}]_{1}=1.$$
\endRoster

\HHH 6.4.4. Vostokov's 
symbol in characteristic 0

\phantom{}\par

Suppose that $\mu _{p^{r}}\le F^*$ and $p>2$. 
Vostokov's symbol
$$ 
V(\;,\;)_{r}\colon K_{m}(F)/p^{r}\times K_{n+1-m}(F)/p^{r}  
\rightarrow K_{n+1}(F)/p^{r} 
\rightarrow  \mu_{p^{r}} 
$$
is defined in section 8.3. For its properties see~8.3. 

\smallskip

 Each pairing defined above is sequentially continuous, so it induces the
pairing of $ K_m^{\tpp}(F)$. 

\HH 6.5. Structure of $K^{{\tpp}}(F)$. I

Denote $VK_{m}^{  {\tpp}}(F)=\bigl\langle \left\{ V_F\right\} \cdot K_{m-1}^   {\tpp}(F)\bigr\rangle $.  Using the tame
symbol and valuation ${\goth v}$
as described in the beginning of 6.4 it is easy to deduce that %
$$
K_{m}(F)\simeq VK_{m}(F)\oplus  \Bbb Z^{a(m)}\oplus ( \Bbb Z%
/(q-1))^{b(m)}$$
with appropriate integer $a(m), b(m)$ (see \cite{FV, Ch. IX, \S2}); 
similar calculations are applicable to $K_{m}^   {\tpp}(F)$.
For example, $\Bbb Z^{a(m)}$ corresponds to 
$\oplus\langle\{t_{j_1},\dots,t_{j_m}\}\rangle$ with $1 \le j_1<\dots<j_m\le n$. 

To study $VK_{m}(F)$ and $VK_{m}^{\tpp}(F)$ the following elementary equality
is quite useful  
$$
\left\{ 1-\alpha ,1-\beta \right\} =\bigl\{ \alpha ,1+\frac{\alpha \beta }{%
1-\alpha }\bigr\} +\bigl\{ 1-\beta ,1+\frac{\alpha \beta }{1-\alpha }%
\bigr\} .
$$
Note that
$\bold v(\alpha\beta/(1-\alpha))=\bold v(\alpha)+\bold v(\beta)$ if
$\bold v(\alpha), \bold v(\beta)>(0,\dots,0)$.

For $\varepsilon ,\eta \in V_F$
one can apply the previous formula to $\left\{ \varepsilon ,\eta
\right\}\in K_{2}^   {\tpp}(F)$ and using the topological convergence
deduce that 
$$\left\{ \varepsilon ,\eta
\right\} 
 =\sum \left\{ \rho _{i},t_{i}\right\} $$ 
with units $\rho
_{i}=\rho _{i}(\varepsilon ,\eta )$  sequentially continuously 
depending on $\varepsilon ,\eta $. 

Therefore $VK_{m}^   {\tpp}(F)$ is
{\it topologically generated} by symbols 
$$
\bigl\{ 1+\theta t_{n}^{i_{n}}\dots t_{1}^{i_{1}},t_{j_{1}}\dots
,t_{j_{m-1}}\bigr\} ,\quad \theta \in \mu _{q-1}.
$$
In particular, $K_{n+2}^   {\tpp}(F)=0$.

 \th Lemma   

$\bigcap_{l \ge 1}lK_{m}(F)\subset \Lambda _{m}(F)$.
\endth

\pf Proof

First, $\bigcap lK_{m}(F)\subset VK _{m}(F)$. 
Let $x\in VK_m(F)$. Write  
$$x=\sum \bigl\{\alpha_J, t_{j_1},\dots,t_{j_{m-1}}\bigr\}
\mod \Lambda_m(F), \quad \alpha_J\in V_F.$$
 Then 
 $$p^{r}x=\sum \bigl\{ \alpha
_{J}^{p^{r}}\bigr\} \cdot \bigl\{ t_{j_1},\dots ,t_{j_{m-1}}\bigr\} +\lambda
_{r}, \quad \lambda _{r}\in \Lambda _{m}(F).$$ 
It remains to apply  property (3) in 6.2.  
\qed\endpf

\HH 6.6. Structure of $K^{{\tpp}}(F)$. II

This subsection 6.6 and the rest of this section  are not required for 
understanding Parshin's class field theory theory
of higher local fields of characteristic $p$ which is 
discussed in section 7.

The next theorem relates the structure of $VK_{m}^{  {\tpp}}(F)$
with the structure of a simpler object.

\th Theorem 1 {{\rm (\cite{F5, Th. 4.6})}} 

 Let  $\chr(K_{n-1})=p$.  The
homomorphism
 $$ 
g\colon \prod_{J}V_{F} \rightarrow VK_{m}(F), \quad 
(\beta _{J}) \mapsto \sum_{J=\{j_1,\dots, j_{m-1}\}}\bigl\{ \beta _{J},t_{j_{1}},\dots
,t_{j_{m-1}}\bigr\} 
 $$ 
induces a homeomorphism between 
 $ \prod V_{F}/g^{-1}(\Lambda _{m}(F))$
endowed with the quotient of the $*$-topology and 
$VK_{m}^{  {\tpp}}(F)$; $g^{-1}(\Lambda _{m}(F))$
is a closed subgroup.

\endth

Since every closed subgroup of $V_F$ is the intersection of
some open subgroups of finite index in $V_F$ (property (7) of 6.2), 
we obtain the following: 

\th Corollary

 \ $\Lambda _{m}(F)=\bigcap $ open  subgroups  of finite index in $%
K_{m}(F)$.
\endth

\rk Remarks

1. If $F$ is of characteristic $p$, 
there is a complete description of the structure of
$K_m^{\tpp}(F)$ in the language of topological generators
and relations due to Parshin (see subsection~7.2).

2. If $\chr(K_{n-1})=0$, then the border homomorphism in  Milnor $K$-theory  
(see for instance \cite{FV, Ch. IX \S2}) 
induces the homomorphism $$VK_m(F)\to VK_m(K_{n-1})\oplus
VK_{m-1}(K_{n-1}).$$ Its kernel is equal to 
the subgroup of  $VK_m(F)$ generated by symbols
$\{u,\dots\}$ with $u$ in the group $1+\Cal M_F$ which is uniquely divisible. 
So  
$$VK_m^{\tpp}(F)\simeq VK_m^{\tpp}(K_{n-1})\oplus VK_{m-1}^{\tpp}(K_{n-1})$$
and one can apply Theorem 1 to describe $VK_m^{\tpp}(F)$. 
\endrk

\pf Proof of Theorem 1

Recall that every symbol
$\{\alpha_1,\dots,\alpha_m\}$ in
$K_m^{\tpp}(F)$ can be written as a convergent sum of
symbols $\{\beta_J, t_{j_1},\dots,t_{j_{m-1}}\}$
with $\beta_J$ sequentially 
continuously 
depending on $\alpha_i$ (subsection~6.5).
Hence there is a sequentially continuous map
$f\colon V_F\times F^{*\oplus m-1}\to \prod_J V_F$
 such that
its composition with $g$ coincides with the restriction of the map 
$\varphi\colon (F^*)^m\to K_m^{\tpp}(F)$ on  $V_F\oplus F^{*\oplus m-1}$.

So the quotient of the $*$-topology  
of $\prod_{J}V_{F}$ is $ \le \lambda _{m}$,
as follows from the definition
of $\lambda_m$. 
Indeed, the sum of  two convergent sequences $x_i,y_i$ in 
$\prod_J V_F/g^{-1}(\Lambda_m(F))$ converges to the sum of their limits. 

Let $U$ be an open subset in   $VK_m(F)$.
Then $g^{-1}(U)$ is open in the $*$-product of the topology $\prod_J V_F$.
Indeed, otherwise for some $J$ 
 there were a sequence 
$\alpha_J^{(i)}\not \in g^{-1}(U)$ which converges to 
$\alpha_J\in g^{-1}(U)$. 
Then the properties of the map  $\varphi$ of 6.3 imply that the sequence
$\varphi(\alpha_J^{(i)})\not\in U$ converges to $\varphi(\alpha_J)\in U$ which contradicts the openness of $U$. 
\qed\endpf

\th Theorem 2 {{\rm(\cite{F5, Th. 4.5})}}

 If $\chr(F)=p$ then $\Lambda _{m}(F)$ is equal to $\bigcap_{l\ge 
1}lK_{m}(F)$ and  is a divisible group.
\endth

\pf Proof 

Bloch--Kato--Gabber's theorem (see subsection~A2 in the appendix to
section~2) shows that the differential symbol
 $$ 
d\colon K_{m}(F)/p \longrightarrow \Omega _{F}^{m}, \qquad 
\left\{ \alpha _{1},\dots ,\alpha _{m}\right\}  \longmapsto \frac{d\alpha
_{1}}{\alpha _{1}}\wedge \dots \wedge \frac{d\alpha _{m}}{\alpha _{m}} $$ 
is injective. \ The topology of $\Omega _{F}^{m}$ induced from $F$ (as
a finite dimensional vector space) is Hausdorff, and $d$ is continuous, so $%
\Lambda _{m}(F)\subset pK_{m}(F)$.

Since $VK_{m}(F)/\Lambda _{m}(F)\simeq \prod \Cal E_{J}$ doesn't have $p$%
-torsion by Theorem~1 in subsection~7.2, $\Lambda _{m}(F)=p\Lambda _{m}(F)$.
\qed\endpf

\th Theorem 3 {{\rm(\cite{F5, Th. 4.7})}} 

If $\chr(F)=0$ then $\Lambda _{m}(F)$ is equal to $\bigcap_{l \ge 1}lK_{m}(F)$ and is a
divisible group. 
If a primitive $l$th root $\zeta_l$  belongs to $F$,  
then  ${}_lK_{m}^   {\tpp}(F)=\left\{ \zeta_l\right\} \cdot K_{m-1}^{\tpp}(F)$.
\endth

\pf Proof

 To show that $p^{r}VK_{m}(F)\supset \Lambda _{m}(F)$ it suffices 
 to check that $p^{r}VK_{m}(F)$ is the intersection of  open neighbourhoods of
$p^{r}VK_{m}(F)$.

We can assume that 
 $\mu_p$ is contained in $F$
applying the standard argument by using $(p,|F(\mu_p):F|)=1$
and $l$-divisibility of $VK_m(F)$ for $l$ prime to $p$.

If $r=1$ then one can use Bloch--Kato's description of
$$
U_{i}K_{m}(F)+pK_{m}(F)/U_{i+1}K_{m}(F)+pK_{m}(F)
$$
in terms of products of quotients of $\Omega _{K_{n-1}}^{j}$
(section~4). 
$\Omega_{K_{n-1}}^j$ and its quotients are finite-dimensional vector spaces 
over $K_{n-1}/K_{n-1}^p$, so the intersection of all neighborhoods of 
zero there 
with respect to the induced from $K_{n-1}$ topology 
is trivial. 
Therefore the injectivity of $d$ implies $\Lambda_m(F)\subset pK_m(F)$. 

Thus, the intersection of open subgroups in $VK_m(F)$ containing
$pVK_m(F)$ is equal to $pVK_m(F)$.

Induction Step.

For a field $F$  consider the pairing
$$(\,\, , \,)_r\colon K_m(F)/p^r\times H^{n+1-m}(F, \mu_{p^r}^{\otimes n-m})\to H^{n+1}(F,\mu_{p^r}^{\otimes n})$$ 
given by the cup product and the map $F^*\to H^1(F,\mu_{p^r})$.
If $\mu_{p^r}\subset F$, then Bloch--Kato's theorem shows that 
$(\,\, ,\,)_r$ can be identified (up to  sign) with  Vostokov's pairing
$V_r(\,,\,)$.  

For $\chi\in H^{n+1-m}(F, \mu_{p^r}^{\otimes n-m})$ put 
$$A_\chi=\{\alpha \in K_m(F): (\alpha,\chi)_r=0\}.$$
One can show \cite{F5, Lemma 4.7} that $A_\chi$ is an open subgroup of  $K_m(F)$.

Let $\alpha$ belong to the intersection of all open subgroups
of $VK_m(F)$ which contain $p^rVK_m(F)$. 
Then $\alpha\in A_\chi$ for every $\chi\in H^{n+1-m}(F,\mu_{p^r}^{\otimes n-m})$. 

Set $L=F(\mu_{p^r})$ and $p^s=|L:F|$.
From the induction hypothesis we deduce that
$\alpha\in p^sVK_m(F)$ and hence $\alpha=N_{L/F}\beta$
for some $\beta \in VK_m(L)$.
Then
$$0=(\alpha,\chi)_{r,F}=(N_{L/F}\beta,\chi)_{r,F}
=(\beta,i_{F/L}\chi)_{r,L}$$
where $i_{F/L}$ is the natural map.
Keeping in mind the identification
between Vostokov's pairing $V_r$ and $(\,\,,\,)_r$ for the field $L$
we see that $\beta$ is annihilated by $i_{F/L}K_{n+1-m}(F)$
with respect to Vostokov's pairing. 
Using explicit calculations with Vostokov's pairing 
 one can directly deduce that 
$$\beta\in (\sigma -1)K_m(L)+p^{r-s} i_{F/L}K_m(F)+p^rK_m(L), $$ 
and therefore $\alpha\in p^rK_m(F)$, as required.

Thus $p^{r}K_{m}(F)= 
\bigcap$ open neighbourhoods of
$p^{r}VK_{m}(F)$.

\smallskip

To prove the second statement we can assume that $l$ is a prime.
If $l\not=p$, 
then 
since $K_m^{\tpp}(F)$ is the direct sum of several cyclic groups and  $VK_m^{\tpp}(F)$ and since $l$-torsion of $K_m^{\tpp}(F)$
is $p$-divisible and $\cap_r p^rVK_m^{\tpp}(F)=\{0\}$,
we deduce 
the result.

Consider the most difficult case of $l=p$. 
Use the exact sequence  $$0\rightarrow \mu _{p^{s}}^{\otimes
n}\rightarrow \mu _{p^{s+1}}^{\otimes n}\rightarrow \mu _{p}^{\otimes
n}\rightarrow 0$$ and the following commutative diagram (see also subsection~4.3.2) 
$$\CD
\mu _{p}\otimes K_{m-1}(F)/p @>>>   K_{m}(F)/p^{s}  @>p>> K_{m}(F)/p^{s+1} \\ 
@VVV @VVV @VVV   \\ 
H^{m-1}(F,\mu _{p}^{\otimes m}) @>>>  H^{m}(F,\mu
_{p^{s}}^{\otimes m})  @>>>  H^{m}(F,\mu _{p^{s+1}}^{\otimes m}).
\endCD
$$
We deduce that $px\in \Lambda _{m}(F)$ implies $px\in \bigcap
p^{r}K_{m}(F)$, so $ x=\left\{ \zeta _{p}\right\} \cdot
a_{r-1}+p^{r-1}b_{r-1}$ for some $a_i\in K_{m-1}^{\tpp}(F)$ and
$b_i\in K_m^{\tpp}(F)$. 

Define $\psi\colon K_{m-1}^{\tpp}(F)\to K_m^{\tpp}(F)$ as $\psi (\alpha )=\left\{ \zeta _{p}\right\} \cdot \alpha $; 
it is a continuous map.  
Put 

\noindent$D_{r}=\psi^{-1}(p^{r}K_{m}^{\tpp
}(F))$. The group $D=\cap D_r$ is the kernel of $\psi$. 
One can show \cite{F5, proof of Th. 4.7} that $\{a_r\}$ is a Cauchy sequence
in the space $K_{m-1}^{\tpp}(F)/D$ which is complete.
Hence  
 there is $y\in
\bigcap \left( a_{r-1}+D_{r-1}\right) $.
Thus,  $x=\left\{ \zeta _{p}\right\}\cdot y$ in 
$K_{m}^   {\tpp}(F)$. 

Divisibility follows.
\qed\endpf

\rk Remarks

1. Compare with  Theorem 8 in 2.5.

2. For more properties of $K_{m}^   {\tpp}(F)$ see \cite{F5}.

3. Zhukov  \cite{Z, \S7--10} gave a description  of $K_{n}^   {\tpp}(F)$ 
in terms of topological generators and relations for some fields $F$
of characteristic zero with small $v_{F}(p)$. 
\endrk

\HH 6.7. The group $K_m(F)/l$

\HHH 6.7.1 

If a prime number $l$ is distinct from $p$, 
then, since $V_F$ is $l$-divisible, we deduce from~6.5 that
$$K_m(F)/l\simeq K_m^{\tpp}(F)/l\simeq (\Bbb Z/l)^{a(m)}\oplus (\Bbb Z/d)^{b(m)}$$
where $d=\text{gcd}(q-1,l)$.

\HHH 6.7.2

The case of $l=p$ is more interesting and difficult.
We use the method described at the beginning of 6.4. 

\smallskip

If $\chr(F)=p$ then
the Artin--Schreier pairing of~6.4.3 for $r=1$
helps one to show that
$
K_{n}^{\tpp}(F)/p$ has the following topological $\Bbb Z/p$-basis: 
$$\bigl\{ 1+\theta t_{n}^{i_{n}}\dots t_{1}^{i_{1}},t_{n},\dots ,\widehat{t_l},\dots ,t_{1}\bigr\}$$
where $p\nmid \text{gcd}(i_1,\dots,i_n)$, \quad $0<(i_{1},\dots ,i_{n})$, \quad 
$l=\minn \left\{ k:p\nmid i_{k}\right\}$ \quad  
and $\theta$ runs over all elements of a fixed basis of $K_0$ over $\Bbb F_p$.

\smallskip

If $\chr(F)=0$, $\zeta _{p}\in F^{* }$,
then using Vostokov's symbol (6.4.4 and 8.3)   one obtains  that 
  $K_{n}^{\tpp}(F)/p$ has the following  topological $\Bbb Z_p$-basis
consisting of elements of two types:
$$
\omega _{* }(j)=\bigl\{ 1+\theta _{* }t_n^{pe_n/(p-1)}\dots
t_1^{pe_1/(p-1)},t_{n},\dots ,\widehat{t_j},\dots ,t_{1}\bigr\}$$
where  $1\le j\le n$, 
$(e_1,\dots,e_n)={\bold v}_{F}(p)$ and $\theta_*\in \mu_{q-1}$ is such that

\noindent$ 1+\theta _{* }t_n^{pe_n/(p-1)}\dots
t_1^{pe_1/(p-1)}$ doesn't belong to $F^{*p}$ 
\par
\noindent and
$$
\bigl\{ 1+\theta t_{n}^{i_{n}}\dots t_{1}^{i_{1}},t_{n},\dots ,\widehat{t_l},\dots ,t_{1}\bigr\}$$
where $p\nmid \text{gcd}(i_1,\dots,i_n)$, \quad $0<(i_{1},\dots ,i_{n})<p(e_1,\dots,e_n)/(p-1) $, 

\noindent  
$l=\minn \left\{ k:p\nmid i_{k}\right\}$,  \quad  
where $\theta$ runs over all elements of a fixed basis of $K_0$ over $\Bbb F_p$.

If $\zeta _{p}\not\in F^{* }$, then pass to the field $F(\zeta_p)$
and then go back, using the fact that the degree of $F(\zeta_p)/F$ is relatively prime to $p$. One deduces that
 $K_{n}^{\tpp}(F)/p$ has the following  topological $\Bbb Z_p$-basis:
$$
\bigl\{ 1+\theta t_{n}^{i_{n}}\dots t_{1}^{i_{1}},t_{n},\dots ,\widehat{t_l},\dots ,t_{1}\bigr\}$$
where $p\nmid \text{gcd}(i_1,\dots,i_n)$, \quad $0<(i_{1},\dots ,i_{n})<p(e_1,\dots,e_n)/(p-1) $, 

\noindent 
$l=\minn \left\{ k:p\nmid i_{k}\right\}$,  \quad  
where $\theta$ runs over all elements of a fixed basis of $K_0$ over $\Bbb F_p$.

\HH 6.8. The norm map on  $K^{\tpp}$-groups

\df Definition

Define the norm map on $K_{n}^  {\tpp}(F)$ as induced 
by $N_{L/F}\colon K_{n}(L)\rightarrow K_{n}(F)$.

Alternatively in characteristic $p$ one can define
the norm map as in 7.4.
\enddf

\HHH 6.8.1

Put $u_{i_{1},\dots ,i_{n}}= U_{i_{1},\dots ,i_{n}}/U_ {i_{1}+1,\dots ,i_{n}}$. 

\th Proposition {{\rm (\cite{F2, Prop. 4.1} and \cite{F3, Prop. 3.1})}}

Let $L/F$ be a cyclic extension of prime degree $l$
such that the extension of the last finite residue fields is trivial. 
Then there is $s$ and a local parameter $t_{s,L}$ of $L$
such that $L=F(t_{s,L})$. 
Let $t_1,\dots, t_n$ be a system of local parameters of $F$, 
then   $t_1,\dots, t_{s,L},\dots, t_n$ is a system of local parameters of $L$.

Let $l=p$. For a generator $\sigma$ of $\Gal(L/F)$ let
$$
\frac{\sigma t_{s,L}}{t_{s,L}}=1+\theta _{0}t_{n}^{r_{n}}\cdots
t_{s,L}^{r_{s}}\cdots t_{1}^{r_{1}}+\cdots 
$$

Then 
\Roster 
\Item{(1)} if $(i_{1},\dots ,i_{n})<(r_{1},\dots ,r_{n})$ then
$$N_{L/F}\colon u_{i_{1},\dots ,i_{n,L}}\to 
u_{pi_{1},\dots ,i_{s},\dots ,pi_{n},F}$$ sends  
$\theta\in K_0$ to $\theta^{p}$;

\Item{(2)} if $(i_{1},\dots ,i_{n})=(r_{1},\dots ,r_{n})$ then 
$$N_{L/F}\colon u_{i_{1},\dots ,i_{n,L}}\to 
u_{pi_{1},\dots ,i_{s},\dots ,pi_{n},F}$$ sends  
$\theta\in K_0$ to $\theta^{p}-\theta\theta_0^{p-1}$;

\Item{(3)} if $(j_{1},\dots ,j_{n})>0$ then 
$$N_{L/F}\colon u_{j_{1}+r_1,\dots ,pj_s+r_s,\dots, j_n+r_n,L}
\to  u_{j_{1}+pr_1,\dots ,j_s+r_s,\dots, j_n+pr_n,F}$$ 
sends 
$\theta\in K_0$ to $-\theta\theta_0^{p-1}$.
\endRoster 
\endth
\pf Proof

Similar to the one-dimensional case \cite{FV, Ch. III \S1}. \qed
\endpf 

\HHH 6.8.2

If $L/F$ is cyclic of prime degree $l$ then 
$$
K_{n}^{\tpp}(L)=\bigl\langle \{ L^{* }\} \cdot
i_{F/L}K_{n-1}^{\tpp}(F)\bigr\rangle
$$
where $i_{F/L}$ is induced by the embedding $F^*\to L^*$.
For instance (we use the notations of section~1), if $f(L|F)=l$ then
$L$ is generated over $F$ by a root of unity of order prime to $p$; 
if $e_i(L|F)=l$, then use the previous proposition.

\th Corollary 1

Let $L/F$ be a cyclic extension of prime degree $l$.
Then $$|K_{n}^{\tpp}(F):N_{L/F}K_{n}^{\tpp}(L)|=l.$$
\endth

If $L/F$ is as in the preceding proposition, then the element 
$$\bigl\{1+\theta _* t_{n}^{pr_{n}}\cdots t_{s,F}^{r_{s}}\cdots t_{1}^{pr_{1}},
t_1,\dots,\widehat{t_s},\dots,t_n\bigr\},$$
where the residue of $\theta_*$ in $K_0$ doesn't belong to the image
of the map  
$$\Cal O_F@>\theta\mapsto \theta^{p}-\theta\theta_0^{p-1}>> \Cal O_F@>>> K_0,$$ 
is  
a generator of $K_{n}^{\tpp}(F)/N_{L/F}K_{n}^{\tpp}(L)$.

If $f(L|F)=1$ and $l\not=p$, then 
$$
\bigl\{\theta_*,t_1,\dots,\widehat{t_s},\dots, t_n \bigr\}
$$
where $\theta_*\in \mu_{q-1}\setminus \mu_{q-1}^l$ 
is a generator of $K_{n}^{\tpp}(F)/N_{L/F}K_{n}^{\tpp}(L)$.

If $f(L|F)=l$, then
$$
\bigl\{t_1,\dots,t_n \bigr\}
$$
is a generator of $K_{n}^{\tpp}(F)/N_{L/F}K_{n}^{\tpp}(L)$.

\th Corollary 2

 $N_{L/F}$(closed  subgroup) is closed
and 
$N_{L/F}^{-1}$(open  subgroup) is open. 
\endth
\pf Proof

Sufficient to show for an extension of prime degree;
then use the previous proposition and Theorem~1 of 6.6.  
\qed
\endpf

\HHH 6.8.3  

\th Theorem 4 {{\rm (\cite{F2,\S4}, \cite{F3,\S3})}}

Let $L/F$ be a cyclic extension of prime degree $l$ with a generator $\sigma $  
then the sequence 
$$
K_{n}^   {\tpp}(F)/l\oplus K_{n}^ {\tpp}(L)/l 
@>i_{F/L}\oplus(1-\sigma )>> 
K_{n}^ {\tpp}(L)/l 
@>{N_{L/F}}>>
K_{n}^ {\tpp}(F)/l
$$
is exact.
\endth

\pf Proof

Use the explicit description of $K_{n}^ {\tpp}/l$ in 6.7.
\qed\endpf

This theorem together with the description of the torsion of $K_n^{\tpp}(F)$
in~6.6
imply:

\th Corollary

 Let $L/F$ be cyclic with a generator $\sigma $ then 
the sequence 
$$K_{n}^   {\tpp}(L) @>{1-\sigma }>>K_{n}^   {\tpp}(L)%
 @>{N_{L/F}}>> K_{n}^   {\tpp}(F)$$ is exact.
\endth

\medskip

\Bib References

\rf{BK} 
S. Bloch and K.  Kato, 
$p$-adic \'etale cohomology,   Inst. Hautes \'Etudes Sci. Publ. Math.
63(1986),  107--152.

\rf{F1} I. Fesenko,
Explicit constructions in local fields,
Thesis, St. Petersburg Univ. 1987.

\rf{F2}
I.  Fesenko, 
Class field theory of multidimensional local fields of 
 characteristic 0, with the residue field of positive characteristic, 
Algebra i Analiz (1991); 
English translation in
St. Petersburg Math. J. 3(1992), 649--678.

\rf{F3}
I.  Fesenko, 
Multidimensional local class field theory II, 
Algebra i Analiz (1991); 
English translation in
St. Petersburg Math. J. 3(1992),   1103--1126. 

\rf{F4} I. Fesenko,
Abelian extensions of complete discrete valuation fields,
Number Theory Paris 1993/94, 
Cambridge Univ. Press,  1996, 
 47--74.

\rf{F5} I. Fesenko, 
Sequential topologies and quotients of the Milnor $K$-groups of higher local
fields,  
preprint, 
 www.maths.nott.ac.uk/personal/ibf/stqk.ps

\rf{FV} I. Fesenko and S. Vostokov, 
Local Fields and Their Extensions,
AMS, Providence RI, 1993.

\rf{K1} K. Kato, A generalization of local class field theory
by using $K$-groups I,
J. Fac. Sci. Univ. Tokyo Sec. IA 26 No.2 (1979), 303--376.

\rf{K2} 
K. Kato, A generalization of local class field theory
by using $K$-groups II,
J. Fac. Sci. Univ. Tokyo Sec. IA 27 No.3 (1980), 603--683.

\rf{K3} K. Kato, 
Existence theorem for higher local class field theory, 
this volume.

\rf{Kn} B. Kahn, {L'anneau de Milnor d'un corps local
a corps residuel parfait}, 
Ann. Inst. Fourier 34(1984), 19--65.

\rf{P1} A. N. Parshin, 
Class fields and algebraic $K$-theory, 
Uspekhi Mat. Nauk  30(1975), 
 253--254; 
English translation in
Russian Math. Surveys.

\rf{P2}
 A. N. Parshin, 
 Local class field theory, 
Trudy Mat. Inst. Steklov.  (1985); 
 English translation in
Proc. Steklov Inst. Math. 1985, issue 3, 157--185. 

\rf{P3} A. N. Parshin, 
Galois cohomology and Brauer group of local fields, 
Trudy Mat. Inst. Steklov. (1990); 
 English translation in
  Proc. Steklov Inst. Math. 1991, issue 4, 191--201.

\rf{S} R. F. Snipes,
Convergence and continuity:
the adequacy of filterbases,
the inadequacy of sequences,
and spaces in which sequences suffice,
Nieuw Arch. voor Wiskunde (3) 23(1975), 8--31.

\rf{V}  S. V. Vostokov, 
   Explicit construction of class field theory for a multidimensional
local field, 
   Izv.  Akad.  Nauk SSSR Ser.  Mat.   (1985) no.2;
English translation in  Math.  USSR Izv.      26(1986),  263--288.

\rf{Z}
I. B. Zhukov, 
Milnor and topological $K$-groups of multidimensional 
complete fields, 
Algebra i analiz (1997);   English translation in  St. Petersburg Math. J.
9(1998), 69--105. 

\endBib

\Coordinates

Department of Mathematics \  
University of Nottingham

Nottingham NG7 2RD England

E-mail: ibf\@maths.nott.ac.uk
\endCoordinates

\vfill
\pagebreak
\end